\documentclass[11pt,a4paper]{article}
\usepackage[utf8]{inputenc}
\usepackage[T1]{fontenc}
\usepackage{mathptmx}
\usepackage{amsmath}
\usepackage{amsfonts}
\usepackage{amssymb}
\usepackage{amsthm}

\usepackage[pdftex]{graphicx}
\usepackage[space]{grffile}
\usepackage[pdftex,linkcolor=black,pdfborder={0 0 0}]{hyperref} 
\usepackage{calc} 

\usepackage{tikz-cd}
\usepackage{tikz}
\usetikzlibrary{matrix,calc}
\usepackage{enumerate}
\usepackage[affil-it]{authblk}
\usepackage{comment}
\usepackage{longtable}
\usepackage{xcolor}
\usepackage{appendix}
\usepackage{booktabs}

\usepackage{esint}
\usepackage{pifont}
\usepackage{csquotes}
\usepackage{mathtools}
\usepackage{fullpage}
\usepackage{euscript}
\usepackage{listings}
\usepackage{esint}
\usepackage{dsfont}

\usepackage[english]{babel}

\usepackage[style=alphabetic, backend=biber]{biblatex}

\renewbibmacro{in:}{}

\AtEveryBibitem{\clearfield{month}}

\usepackage[a4paper, lmargin=0.1666\paperwidth, rmargin=0.1666\paperwidth, tmargin=0.1111\paperheight, bmargin=0.1111\paperheight]{geometry} 

\usepackage[all]{nowidow} 
\usepackage[protrusion=true,expansion=true]{microtype} 

\hypersetup{
pdfsubject = {},
pdftitle = {A note on Einstein 4-manifolds with symmetries},
pdfauthor = {Patrick Donovan}
}

\numberwithin{equation}{section}

\newtheorem{theorem}{Theorem}[section]
\newtheorem{lemma}[theorem]{Lemma}

\newtheorem{proposition}[theorem]{Proposition}

\newtheorem{corollary}[theorem]{Corollary}

\theoremstyle{definition}
\newtheorem{definition}[theorem]{Definition}

\theoremstyle{remark}
\newtheorem{remark}[theorem]{Remark}

\theoremstyle{plain}
\newtheorem{thmint}{Theorem}

\newcommand{\acknowledgements}[1]{
    \begin{center}
        \bfseries Acknowledgements
    \end{center}
    {\small #1} 
}

\title{A note on cohomogeneity one Einstein metrics on closed 4-manifolds}
\author{Patrick Donovan}

\affil{\small Mathematical Institute, University of Oxford\\ \emph{E-mail address}: {patrick.donovan@balliol.ox.ac.uk}}

\date{November 2026}

\begin{document}

\maketitle

\begin{abstract}
Broder recently classified the closed Einstein $4$-manifolds whose isometry groups have dimension at least three \cite{Bro26}. The proof reduces the problem to the study of cohomogeneity one Einstein metrics invariant under $\mathrm{SU}(2)$, however a few cases (following \cite[Table 4]{Don25}) are not treated. We reconcile the missing cases using techniques from \cite{Bro26} with minimal modifications.
\end{abstract}

\section{Introduction}

Broder \cite[Theorem 1.1]{Bro26} states that if $(M^4,g)$ is a connected closed Einstein 4-manifold with $\dim\mathrm{Isom}(M,g)\geq3$, then $g$ is flat or its universal cover is homothetic to one of
\begin{equation}\label{eq:list}
(\mathbb{S}^4,g_{\mathrm{round}}),\qquad(\mathbb{CP}^2,g_{\mathrm{FS}}),\qquad(\mathbb{S}^2\times \mathbb{S}^2,g_{\mathrm{round}}\oplus g_{\mathrm{round}}),\qquad(\mathbb{CP}^2\#\overline{\mathbb{CP}^2},g_{\mathrm{Page}}).
\end{equation}
First, \cite[Theorem 3]{Der83} (see also \cite{BB82}) implies that only cohomogeneity one actions with 3-dimensional isometry group need to be considered. These actions are by $\mathrm{SU}(2)$ and are known, see \cite{Par86,Hoe10}. In this setting, every invariant Einstein metric is diagonal (c.f. \cite[Theorem 3.1]{Don25}); along a normal geodesic it can be uniquely given in the form $dt^2+\sum_if_i^2\sigma_i^2$. The topologies and smoothness conditions depend on which $f_i$ vanishes at each singular orbit, and the possibilities are listed in \cite[Table 4]{Don25} (Table \ref{tab:smooth} below).

Let $M_n$ be the space whose group diagram (see \S2) is $\mathbb{Z}_n\subset\{\mathrm{U}(1),\mathrm{U}(1)\}\subset\mathrm{SU}(2)$, ($M_{\mathrm{odd}}\cong \mathbb{CP}^2\#\overline{\mathbb{CP}^2}$ and $M_{\mathrm{even}}\cong \mathbb{S}^2\times \mathbb{S}^2$). For $n \in \{1,2\}$, $M_n$ appears twice in Table \ref{tab:smooth} as there are two configurations: either the same warping function vanishes at both ends (row (R5)), or different functions do (row (R6)). In \cite[proof of Theorem 5.1]{Bro26} row (R6) on $\mathbb{CP}^2\#\overline{\mathbb{CP}^2}$ is reduced to row (R5) by changing the group diagram by conjugating one singular orbit by an element of $N(H)_0$, where $H$ is a principal isotropy. In \cite[\S6]{Bro26}, on $\mathbb{S}^2\times \mathbb{S}^2$, only row (R6) is treated. This does not suffice since rows (R5) and (R6) genuinely must be treated separately; see \S \ref{sec:A}. In this note, we close the remaining cases.

\begin{thmint}\label{thm:A}
There is no $\mathrm{SU}(2)$-invariant Einstein metric on $\mathbb{S}^2\times \mathbb{S}^2$ in row (R5), and there is no $\mathrm{SU}(2)$-invariant Einstein metric on $\mathbb{CP}^2\#\overline{\mathbb{CP}^2}$ in row (R6). In particular, up to rescaling, the Page metric and the sum of round metrics are the only Einstein metrics on $M_n$ invariant under the defined $\mathrm{SU}(2)$-action.
\end{thmint}

We only need to consider $n \in \{1,2,4\}$, as otherwise the metric is $\mathrm{U}(2)$-invariant and thus classified in \cite{BB82}. One should note that only techniques from \cite{Bro26} are required to prove Theorem \ref{thm:A}.

We introduce the required notation and preliminary results in \S\ref{sec:prelim}. In \S\ref{sec:A}, we show that the group diagram read off along a normal geodesic is determined by the metric up to interchanging $K_\pm$ and simultaneous conjugation; see Proposition \ref{prop:classes}. Finally, we prove Theorem \ref{thm:A} in \S\ref{sec:S2S2} and \S\ref{sec:CP2}, the former dealing with $\mathbb{S}^2\times\mathbb{S}^2$ and the latter with $\mathbb{CP}^2\#\overline{\mathbb{CP}^2}$.

\acknowledgements{I discussed these issues with Kyle Broder 
after \cite{Bro26} appeared on arXiv. 
While working on a related project, I realised that the remaining cases could be resolved quite simply using the techniques of \cite{Bro26}, which led to this short note. 
I would also like to thank Qiu Shi Wang for interesting conversations related to this topic.}

\section{Preliminaries}\label{sec:prelim}

One familiar with \cite{Bro26}, or cohomogeneity one actions more generally, should skip this section.

Suppose $(M^4,g)$ is a closed simply-connected 4-manifold on which $G=\mathrm{SU}(2)$ acts isometrically and by cohomogeneity one. The orbit space is homeomorphic to a closed interval. Let $c:[0,T]\to M$ be a unit speed geodesic between the singular obits $S_- = G \cdot c(0)$ and $S_+ = G\cdot c(T)$ that intersects each orbit orthogonally and exactly once; we call such geodesics normal. The principal isotropy is $H=G_{c(t)}$ for $t \in (0,T)$, and the singular isotropies are $K_-=G_{c(0)}$, $K_+=G_{c(T)}$. The diagram $H\subset\{K_-,K_+\}\subset G$ determines $M$. The diagrams obtained by interchanging $K_\pm$, simultaneous conjugation of $H,K_{\pm}$ by an element of $G$, and replacing $K_+$ by $aK_+a^{-1}$ with $a\in N(H)_0$ (c.f. \cite{GWZ08}) also determine $M$, up to equivariant diffeomorphism. 

We will use notation introduced in \cite[\S2]{Bro26}. Fix an $\mathrm{Ad}$-invariant inner product $Q$ on $\mathfrak{su}(2)$ and a $Q$-orthonormal basis $\{X_i\}_{i=1}^3$ with $[X_i,X_{i+1}]=2X_{i+2}$, with $\{\sigma_i\}$ a dual coframe of the action fields $\{X_i^*\}$. A diagonal metric is given by $g=dt^2+\sum_if_i(t)^2\sigma_i^2$ on the principal part. It is known that any Einstein metric on $g$ can be written in a diagonal form, and that the functions $f_i$ are uniquely determined up to permutation (c.f. \cite[Theorem 3.1]{Don25}).

For $\{i,j,k\} = \{1,2,3\}$ define
\[
L_i=\frac{f_i'}{f_i},\quad R_i=\frac{f_i}{f_jf_k},\quad S=\sum_iL_i,\quad A_i=L_i+R_j+R_k-R_i,
\]
\[
a_i=2R_iA_i-A_jA_k,\qquad u_{ij}=\log(f_i/f_j).
\]
The Einstein equations, $\mathrm{Ric}(g)=\lambda g$, imply the following equations \cite{Bro26}:
\begin{align}
L_i'&=-SL_i+2R_i^2-2(R_j-R_k)^2-\lambda,\qquad R_i'=R_i(L_i-L_j-L_k),\label{eq:E}\\
A_i'&=(R_j+R_k-3R_i-A_i)A_i+A_jA_k,\qquad a_i'=A_j(a_k-a_i)+A_k(a_j-a_i),\label{eq:Aa}\\
u_{ij}'' + Su_{ij}' &= 4(R_i-R_j)(R_i+R_j-R_k) = 4\frac{(f_i^2-f_j^2)(f_i^2+f_j^2-f_k^2)}{f_i^2f_j^2f_k^2}.\label{eq:u}
\end{align}

Only one smoothness condition is needed in this work: a diagonal metric $g$ on $M_n$ is smooth about $S_-$ if and only if there are smooth functions $\phi_i:(-\varepsilon,\varepsilon) \to \mathbb{R}$, such that \cite{VZ22}
\begin{equation}\label{eq:smooth}
f_1^2=n^2t^2+t^4\phi_1(t^2),\qquad f_{2,3}^2=\phi_2(t^2)\pm t^{4/n}\phi_3(t^2),
\end{equation}
where $\phi_2$ is positive and $\phi_3$ necessarily vanishes unless $n\in\{1,2,4\}$. The smoothness condition at $S_+$ is given similarly; substitute $t\mapsto T-t$ and permute indices if required. Table \ref{tab:smooth} lists all possible smoothness conditions:

\begin{table}[h]
\centering
\renewcommand{\arraystretch}{1.2}
\begin{tabular}{@{}lllll@{}}
\toprule
& Action & $t\to0$ & $t\to T$ & Einstein metrics\\
\midrule
(R1) & $\mathrm{SO}(3)$ on $S^4$ & \eqref{eq:smooth}, $n=4$ & \eqref{eq:smooth}, $n=4$, $f_1\leftrightarrow f_3$ & round \cite[Thm. 3.2]{Bro26}\\
(R2) & $\mathrm{SO}(3)$ on $\mathbb{CP}^2$ & \eqref{eq:smooth}, $n=4$ & \eqref{eq:smooth}, $n=2$, $f_1\leftrightarrow f_3$ & Fubini--Study \cite[Thm. 4.2]{Bro26}\\
(R3) & $\mathrm{SU}(2)$ on $S^4$ & fixed point & fixed point & round \cite[Thm. 3.1]{Bro26}\\
(R4) & $\mathrm{SU}(2)$ on $\mathbb{CP}^2$ & fixed point & \eqref{eq:smooth}, $n=1$ & Fubini--Study \cite[Thm. 4.1]{Bro26}\\
(R5) & $M_n$, $n$ odd & \eqref{eq:smooth} & \eqref{eq:smooth} & Page \cite[Thm. 5.1]{Bro26}\\
(R5) & $M_n$, $n$ even & \eqref{eq:smooth} & \eqref{eq:smooth} & none, Theorem \ref{thm:S2S2}\\
(R6) & $M_2$ & \eqref{eq:smooth} & \eqref{eq:smooth}, $f_1\leftrightarrow f_3$ & product \cite[Thm. 6.1]{Bro26} \\
(R6) & $M_1$ & \eqref{eq:smooth} & \eqref{eq:smooth}, $f_1\leftrightarrow f_3$ & none, Theorem \ref{thm:CP2}\\
\bottomrule
\end{tabular}
\caption{The configurations of \cite[Table 4]{Don25} and where they are treated.}
\label{tab:smooth}
\end{table}

\section{Equivariant diffeomorphisms and the smoothness conditions}\label{sec:A}

This section is a clarification of \cite[Remark 2.9]{Don25}.

\begin{definition}
A $G$-invariant metric $g$ is \emph{compatible} with the group diagram $H\subset\{K_-,K_+\}\subset G$ if some normal geodesic $c$ of $g$ realises it, that is,
\[
G_{c(0)}=K_-,\qquad G_{c(T)}=K_+,\qquad G_{c(t)}=H\ \text{ for all } t\in(0,T).
\]
\end{definition}

There may exist metrics compatible with some group diagram that are not compatible with others representing the same action. In particular, a metric $g$ being compatible with a given group diagram is only preserved by interchanging $K_\pm$ and simultaneous conjugation. That is, the diagram obtained by replacing $K_+$ by $aK_+a^{-1}$ with $a\in N(H)_0$ (c.f. \cite{GWZ08}) has the same space of invariant metrics, but need not have the same space of compatible metrics.

\begin{proposition}\label{prop:classes}
Let $g$ be compatible with $H\subset\{K_-,K_+\}\subset G$.
\begin{enumerate}[(a)]
\item $g$ is compatible with another group diagram $H'\subset\{K'_-,K'_+\}\subset G$ if
and only if there exists $b\in G$ such that $H'=bHb^{-1}$ and
\[
(K'_-,K'_+)=(bK_-b^{-1},\,bK_+b^{-1})\qquad\text{or}\qquad(K'_-,K'_+)=(bK_+b^{-1},\,bK_-b^{-1}).
\]
That is, the new diagram is obtained by conjugating all three isotropy subgroups by the same
element, possibly after interchanging $K_-$ and $K_+$.
\item If $\varphi\colon M'\to M$ is a $G$-equivariant diffeomorphism, then $\varphi^*g$
is compatible with the same group diagrams as $g$.
\end{enumerate}
\end{proposition}

\begin{proof}
Let $c$ be a normal geodesic of $g$ realising $H\subset\{K_-,K_+\}\subset G$.

(a): If the new diagram is of this form, then $b\cdot c$ (or $t \mapsto b \cdot c(T-t)$ if $K_\pm$ are interchanged) is a normal geodesic of $g$ realising the new diagram. Conversely, suppose $g$ is compatible with a group diagram that is realised by $\tilde c$ a normal geodesic of $g$. Reversing $\tilde c$ if necessary, we may assume that $\tilde c(0)\in S_-$, say $\tilde c(0)=h\cdot c(0)$ with $h\in G$. Then $h^{-1}\tilde c$
is a normal geodesic starting at $c(0)$. By the slice theorem, $K_- = G_{c(0)}$ acts transitively on $\partial B (T_{c(0)}S_-^\perp)$, i.e., the space of unit vectors perpendicular to $S_-$ (c.f. \cite[proof of Lemma 2.6]{Don25}). Thus, there is a $k\in K_-$ that maps the initial velocity of $c$ to that of $h^{-1} \tilde{c}$, i.e., $c'(0) = (k^{-1}h^{-1}\tilde c)'(0)$. By uniqueness of geodesics, $\tilde c=b\cdot c$ with $b=hk$, and
$G_{\tilde c(t)}=b\,G_{c(t)}\,b^{-1}$ for all $t\in[0,T]$.

(b): The map $\varphi$ is an isometry from $(M',\varphi^*g)$ to $(M,g)$ mapping orbits to
orbits, so $\varphi^{-1}\circ c$ is a normal geodesic of $\varphi^*g$. By equivariance,
$G_{\varphi(p)}=G_p$, so $\varphi^{-1}\circ c$ realises the same group diagram as $c$.
\end{proof}

Suppose $G$ acts via cohomogeneity one on $(M,g)$ and let $c$ be a normal geodesic. Let $\mathfrak{g}$ be the Lie algebra of $G$ and let $\mathfrak{h}$ be the Lie algebra of $H$, a principal isotropy. Let $\mathfrak{n}$ be an $Ad_H$-invariant subspace of $\mathfrak{g}$ such that $\mathfrak{n} \oplus \mathfrak{h} = \mathfrak{g}$. Identifying on the principal part of $M$ with  $(G/H)\times(0,T)$ via $(h,t)\mapsto h\cdot c(t)$, the metric takes the form $g=dt^2+g_t$, where
$g_t(X,Y)=g(X^*,Y^*)|_{c(t)}$ is an inner product on $\mathfrak{n}$ for each $t \in (0,T)$. The following well known and simple result specifies to which group diagram a specific metric belongs along $c$.

\begin{lemma}\label{lem:kernel}
$(g_t)_{t \in (0,T)}$ extend continuously to positive semi-definite forms $g_0$ and $g_T$ satisfying
\[
\mathfrak{k}_- = \ker g_0 \oplus \mathfrak{h} \qquad\text{and}\qquad \mathfrak{k}_+=\ker g_T \oplus \mathfrak{h},
\]
where $\mathfrak{k}_\pm$ are the Lie algebras of $K_-=G_{c(0)}$ and $K_+=G_{c(T)}$.
\end{lemma}

\begin{proof}
We consider $t=0$; the argument at $t=T$ is the same. $g_t$ extends continuously to a positive semi-definite form $g_0$ since for any $X \in \mathfrak{n}$, $g_t(X,X)=|X^*(c(t))|^2 \underset{t \to 0}{\to}|X^*(c(0))|^2$. An element $X \in \mathfrak{n}$ also lies in $\mathfrak{k}_- $ precisely when 
$X^*(c(0))=0$, i.e., $X\in\ker g_0$. That is, $\mathfrak{n} \cap \mathfrak{k}_- = \ker g_0$ and the result now follows since $\mathfrak{h} \subset \mathfrak{k}_-$.
\end{proof}

We now specialise to the case of cohomogeneity one actions on 4 manifolds. In particular, due to the Einstein condition, we may assume that the metric is diagonal. Recall that we only consider the case where $G= \mathrm{SU}(2)$, $K_\pm$ are $\mathrm{U}(1)$ subgroups of $\mathrm{SU}(2)$, and $H = \mathbb{Z}_n$.

\begin{corollary}\label{cor:diag}
Let $g=dt^2+\sum_if_i(t)^2\sigma_i^2$ be diagonal with respect to a normal geodesic $c$ and a basis $X_1,X_2,X_3$ of $\mathfrak{su}(2)$ as in \S\ref{sec:prelim}.
If $f_i$ vanishes at $t=0$ and $f_j$ vanishes at $t=T$, then $\mathfrak{k}_-=\mathrm{span}(X_i)$ and $\mathfrak{k}_+=\mathrm{span}(X_j)$.
\end{corollary}

\begin{proof}
By Lemma \ref{lem:kernel}, $\mathfrak{k}_- = \ker g_0$ is spanned by the $X_k$ for which $f_k(0)=0$.
\end{proof}

We say a $G$-invariant metric $g = dt^2 +g_t$ on $M$ is diagonalisable if there exists a single choice of basis of $\mathfrak{n}$ as in \S 2 such that $g_t$ is diagonal for all $t \in (0,T)$.
\begin{corollary}\label{thm:A}
We say $g \in M_{n,s}^G$ if $g$ is a diagonalisable $G$-invariant metric on $M_n$ such that for some choice of geodesic and some choice of basis of $\mathfrak{su}(2)$, the same can function vanishes at both ends. Similarly, we say $g \in M_{n,d}^G$ if different functions can vanish. $M_{n,d}^G$ and $M_{n,s}^G$ are disjoint and invariant under pull-back by $G$-equivariant diffeomorphisms $M_n \to M_n$.

The Page metric belongs to $M_{1,s}^G$ (R5) and the Riemannian sum of round metrics belongs to $M_{2,d}^G$ (R6).
\end{corollary}

\begin{proof}
Note that the row to which each metric belongs follows from Corollary \ref{cor:diag}. The latter also implies that $\mathfrak{k}_- = \mathfrak{k}_+$ for some group diagram if and only if $g \in M_{n,s}^G$. Since the condition $\mathfrak{k}_- = \mathfrak{k}_+$ is preserved under conjugation, the result follows immediately from Proposition \ref{prop:classes}.
\end{proof}

\begin{remark}
    In general, to study the space of all cohomogeneity one metrics on $M$ of the form $dt^2 + g_t$, Proposition \ref{prop:classes} shows that one must consider many different group diagrams. In particular, given an action with some group diagram $H \subset \{K_-,K_+\} \subset G$, one must check the space of metrics compatible with the diagram $H \subset \{K_-,aK_+a^{-1}\} \subset G$ for  $a \in N(H)_0$. It suffices to let $a$ range over a set of representatives of the double cosets of
    \[L_-\backslash (L_- \cdot N(H)_0 \cdot L_+ )/L_+, \qquad \text{where} \; L_{\pm} := N(H) \cap N(K_{\pm}).\]
    This formulation is essentially given in \cite{Neu68}. If $K_-$ and $K_+$ are not conjugate to each other by an element of $N(H)$, then each double coset represents a different class of compatible invariant metrics. If $K_-$ and $K_+$ are conjugate to each other by an element $c \in N(H)$, then one could also account for reversing the normal geodesic. That is, the diagrams $H \subset \{K,cKc^{-1}\} \subset G$ and $H \subset \{K,bKb^{-1}\} \subset G$ for $b \in N(H)_0c$ or $b \in N(H)_0 c^{-1}$ have the same space of compatible metrics if and only if
    \[b \in N(K) \cap N(H) \cdot c \cdot N(K) \cap N(H) \]
    or 
    \[b \in N(K) \cap N(H) \cdot c^{-1} \cdot N(K) \cap N(H) \]
    respectively.
\end{remark}

\section{$\mathbb{S}^2\times \mathbb{S}^2$ in row (R5)}\label{sec:S2S2}

\begin{theorem}\label{thm:S2S2}
For $n$ even, $M_n$ carries no $\mathrm{SU}(2)$-invariant Einstein metric in row (R5).
\end{theorem}

This follows using the exact same treatment as in Theorem 5.1 (see also Theorem 3.1 and Theorem 4.1) of \cite{Bro26} (which we repeat below for the convenience of the reader).

\begin{proof}
Suppose $f_1$ vanishes at both ends. By \eqref{eq:smooth}, as $t \to 0$ and $t \to T$, $f_2/f_3 \to 1$ and $f_1/f_3\to0$. Suppose some ratio $f_i/f_3$ with $i \in \{1,2\}$ exceeds $1$ somewhere, and hence is attained at an interior point $t_0$. Let $j \in \{1,2\}$, $j \neq i$ be the remaining index. At $t_0$, we have that $f_j(t_0) \leq f_i(t_0)$, $0<f_3(t_0) < f_i(t_0)$ and $u_{i3}=\log(f_i/f_3)$ reaches an interior maximum, but \eqref{eq:u}
gives
\[u_{i3}''(t_0)=4\,\frac{(f_i^2-f_3^2)(f_i^2+f_3^2-f_j^2)}{f_i^2f_j^2f_3^2}\Big|_{t_0}>0,\]
a contradiction. Thus the ratio $f_2/f_3$ does not exceed $1$, i.e., $f_2 \leq f_3$. Similarly, $f_3 \leq f_2$, so $f_2\equiv f_3$. This implies that the metric is $\mathrm{U}(2)$-invariant and we can appeal to the classification of $\mathrm{U}(2)$-invariant compact Einstein manifolds \cite{BB82}.
\end{proof}

\section{$\mathbb{CP}^2\#\overline{\mathbb{CP}^2}$ in row (R6)}\label{sec:CP2}

We now rule out Einstein metrics in row (R6). This is achieved by using the conditions at the endpoints to contradict the forward invariant cone $\{a_1-a_2, a_2 - a_3 > 0\}$. To show that this cone is forward invariant, we use the maximum principal from the previous section and the positivity of $A_i$. This follows along the lines of the proof of Theorem 3.2 of \cite{Bro26}.

Let $g = dt^2 + \sum_{i=1}^3f_i(t)^2 \sigma_i^2$ be an $\mathrm{SU}(2)$-invariant Einstein metric on $M_1$ in row (R6), i.e., with $f_1$ vanishing at $t=0$ and $f_3$ at $t=T$. Recall that a compact Einstein manifold with a non-abelian identity component of the isometry group implies the Einstein constant is positive. As in \cite{Don25}, define the constants $\beta_-=1/f_2(0)^2=1/f_3(0)^2$ and $\beta_+=1/f_2(T)^2 = 1/f_1(T)^2$, and also parametrise the reverse direction by $s$; $s=T-t$. Solving \eqref{eq:E} order by order about $t=0$ using \eqref{eq:smooth} gives,
\[L_1 = \frac{1}{t} + O(t), \; L_2,L_3 = \frac{4\beta_- - \lambda}{2} t + O(t^3), R_1 = \beta_- t + O(t^3), R_2,R_3 = \frac{1}{t}+\frac{2\beta_-}{3}+ O(t^3),\]
and hence around $s=0$, 
\[L_1,L_2 = \frac{\lambda - 4\beta_+}{2} s + O(s^3), L_3 = -\frac{1}{s} + O(s), \; R_1,R_2 = \frac{1}{s} + \frac{2\beta_+}{3}+ O(s^3),\; R_3 = \beta_+ s + O(s^3).\]
In particular, about $t=0$, 
\[A_1 = \frac{3}{t} + O(t), \; A_2, A_3 = \frac{6\beta_- - \lambda}{2}t + O(t^3), \; a_1 = 6 \beta_- +O(t^2), \; a_2, a_3 = \frac{\lambda-6 \beta_-}{2} + O(t^2),\]
and about $s=0$,
\[A_1,A_2 = \frac{\lambda - 2\beta_+}{2}s + O(s^3), \; A_{3} = \frac{1}{s} + O(s),\; a_1,a_2 = \frac{\lambda-2 \beta_+}{2} + O(s^2), \; a_3 = 2 \beta_+ +O(s^2).\]

\begin{lemma}\label{lem:f_bounds}
    On $(0,T)$, $f_2 > f_1,f_3$. In particular $R_2 > R_1, R_3$ on $(0,T)$.
\end{lemma}

\begin{proof}
    By applying a similar argument to the one in Theorem \ref{thm:S2S2}, $f_3/f_2,f_1/f_2 \leq 1$. If, at some $t_0 \in (0,T)$, $f_i(t_0) = f_2(t_0)$ for $i \neq 2$, then $f_2 \geq f_i$ implies that $f_i'(t_0) = f_2(t_0)$, but then uniqueness implies that $f_i \equiv f_2$, which contradicts the initial conditions in (R6).
\end{proof}

\begin{lemma} \label{lem:lambda_bounds}
    $ 0 < \lambda < 6\beta_-,6 \beta_+$, and on $(0,T)$, $A_i >0$.
\end{lemma}

\begin{proof}
Equation \eqref{eq:Aa} implies that for $\epsilon_i \in \{\pm1\}$ such that $\epsilon_1 \epsilon_2 \epsilon_3 = 1$ the cone $\{\epsilon_i A_i >0\}$ is forward invariant. By the initial and final conditions, we have that $A_1(t)>0$ for small $t$ and $A_3(s) > 0$ for small $s$. In particular, if $A_2'(0) \neq 0$, then, since $A_2$ and $A_3$ share the same sign for sufficiently small $t>0$, they must both be positive. If $A_2'(0) = A_3'(0 ) = 0$, then they both contain at the lowest third order terms. However, using the initial conditions and \eqref{eq:Aa},
\[\left( \begin{matrix}
    A_2 \\
    A_3
\end{matrix}\right)' = \frac{1}{t}Q \left( \begin{matrix}
    A_2 \\
    A_3
\end{matrix}\right)+ B(t) \left( \begin{matrix}
    A_2 \\
    A_3
\end{matrix}\right)\]
for some bounded matrix $B(t)$ and $Q = \left( \begin{matrix}
    -2 &3 \\
    3 & -2
\end{matrix}\right)$. The eigenvalues of $Q$ are $-5$ and $1$, so by Lemma 2.5 of \cite{Bro26}, $A_2 \equiv 0 \equiv A_3$, for small $t>0$ and hence on $(0,T)$ by uniqueness, which again contradicts the behaviour of $A_3$ near $T$. 
By symmetry of the Einstein equations, $\beta_+$ must also satisfy this bound.
\end{proof}

Note that \eqref{eq:Aa} implies that
\begin{equation} \label{eq:A}
    (A_i-A_j)' = (R_j + R_k - 3 R_i -A_i - A_j - A_k) (A_i - A_j) + 4(R_j-R_i)A_j.
\end{equation}

\begin{lemma} \label{lem:A}
    On $(0,T)$, $A_2 < A_1,A_3$.
\end{lemma}

\begin{proof}
    Recall that on $(0,T)$, $R_2 > R_1,R_3$ and $A_i >0$ by Lemmas \ref{lem:f_bounds} and \ref{lem:lambda_bounds}. By \eqref{eq:A}, $A_i - A_j > 0$ is preserved, so it suffices to show the lemma holds for sufficiently small $t>0$. This is trivial for $A_1$ by the initial conditions. For $A_3$, we integrate \eqref{eq:A} to see that for positive $t_0$, $(A_3 - A_2)(t_0) = e^{\int c} \int_0^{t_0} 4e^{-\int c}(R_2-R_3)A_2 \; dt >0$ where $c = R_1+R_2-3R_3-A_1-A_2-A_3$.
\end{proof}

\begin{theorem}\label{thm:CP2}
There is no $\mathrm{SU}(2)$-invariant Einstein metric on $\mathbb{CP}^2\#\overline{\mathbb{CP}^2}$ in row (R6).
\end{theorem}

\begin{proof}
By \eqref{eq:Aa},
\[(a_1-a_2)' = -(2A_3+A_2)(a_1-a_2) + (A_1-A_2)(a_2-a_3)\]
\begin{equation} \label{eq:a}
    (a_2-a_3)' = -(2A_1+A_2)(a_2-a_3) + (A_3-A_2)(a_1-a_2).
\end{equation}
Lemma \ref{lem:A} implies that the cone $a_1-a_2,a_2-a_3 > 0$ is forward invariant on $(0,T)$. The initial conditions and Lemma \ref{lem:lambda_bounds} implies that $(a_1-a_2)(0) = \frac{18\beta_- - \lambda}{2} > 0$ and $(a_2-a_3)(T) = \frac{\lambda - 6\beta_+}{2}<0$. Thus, for a contradiction, it suffices to show that $a_2-a_3$ is positive at some time. As before, we integrate \eqref{eq:a} to see that for positive $t_0$:
\[(a_2 - a_3)(t_0) = e^{-\int 2A_1+A_2} \int_0^{t_0} e^{\int 2A_1+A_2}(A_3-A_2)(a_1-a_2) \; dt >0.\]
\end{proof}

\printbibliography

\end{document}